\documentclass[12pt]{amsart}
\usepackage[margin=1in]{geometry}
\usepackage[utf8]{inputenc}
\usepackage{amsmath}
\usepackage{amssymb}
\usepackage{esint}

\usepackage[pdftitle={HoloH2C},
  pdfauthor={Samuel Bronstein},
  pdffitwindow=true,
  breaklinks=true,
  colorlinks=true,
  urlcolor=blue,
  linkcolor=red,
  citecolor=red,
  anchorcolor=red]{hyperref}
\usepackage{titleps}
\usepackage{xcolor}
\definecolor{links}{RGB}{70,0,255}
\hypersetup{citecolor=blue,linkcolor=links}

\usepackage{subfiles}

\numberwithin{equation}{section}
\numberwithin{table}{section}
\numberwithin{figure}{section}

\newcommand{\C}{\mathbb C}

\newcommand{\R}{\mathbb R}

\newcommand{\Z}{\mathbb Z}
\renewcommand{\H}{\mathbb H}

\newcommand{\Hom}{\text{Hom}}
\newcommand{\Tr}{\text{Tr}}

\newcommand{\calE}{\mathcal E}
\newcommand{\calF}{\mathcal F}

\newcommand{\calO}{\mathcal O}

\newcommand{\calU}{\mathcal U}
\newcommand{\calV}{\mathcal V}

\newcommand{\dbar}{\overline{\partial}}
\newcommand{\eps}{\varepsilon}
\newcommand{\End}{\text{End}}
\newcommand{\Vol}{\text{Vol}}

\newcommand{\II}{\mathrm{I\!I}}

\newcommand{\PO}{\mathrm{PO}}
\newcommand{\PU}{\mathrm{PU}}

\newcommand{\SO}{\mathrm{SO}}
\newcommand{\SU}{\mathrm{SU}}

\newcommand{\PSL}{\mathrm{PSL}}
\newcommand{\one}{\mathbf{1}}
\newcommand{\Isom}{\mathrm{Isom}}

\theoremstyle{plain}
\newtheorem{theorem}{Theorem}[section]
\newtheorem*{theorem*}{Theorem}
\newtheorem{lemma}[theorem]{Lemma}
\newtheorem{proposition}[theorem]{Proposition}

\newtheorem{corollary}[theorem]{Corollary}
\newtheorem*{corollary*}{Corollary}

\theoremstyle{definition}
\newtheorem{definition}[theorem]{Definition}

\newtheorem{notation}[theorem]{Notation}
\newtheorem{question}[theorem]{Question}

\theoremstyle{remark}
\newtheorem*{remark}{Remark}

\title{Almost-fuchsian representations in $\PU(2,1)$}
\author{samuel bronstein}
\date{2023}

\begin{document}
\begin{abstract}
	In this paper, we study nonmaximal representations of surface groups in $\PU(2,1)$.
	In genus large enough, we show the existence of convex-cocompact representations
	of non-maximal Toledo invariant admitting a unique equivariant minimal surface,
	which is holomorphic and almost totally geodesic.
	These examples can be obtained for any Toledo invariant of the form $2-2g+\frac{2}{3}d$,
	provided $g$ is large compared to $d$. When $d$ is not divisible by $3$,
	this yields examples of convex-cocompact representations in $\PU(2,1)$ which do not lift
	to $\SU(2,1)$.
\end{abstract}
\maketitle
\tableofcontents
\section{Introduction}
Let $\Sigma$ be a closed oriented surface of genus $g$ at least $2$ endowed with a hyperbolic
metric. Let $G=\PU(2,1)$ the isometry group of the complex hyperbolic plane.
Xia \cite{Xia00} counted the connected components of the character variety of $\pi_1\Sigma$ into $\PU(2,1)$.
He proved that there are $6g-5$ connected components, indexed by the Toledo invariant $\tau(\rho)$ which belongs to $\frac{2}{3}\Z$
and satisfies a Milnor--Wood type inequality:
\begin{equation}
2-2g\leq \tau(\rho)\leq 2g-2\,.
\end{equation}
When $|\tau(\rho)|=2g-2$, the representation is maximal, and preserves a totally geodesic copy of $\H^1_\C$ in $\H^2_\C$, cf Toledo~\cite{Tol79}.
This has later been generalized to representations in any Hermitian group, see Burger--Iozzi--Wienhard and Burger--Iozzi--Labourie--Wienhard~\cite{BIW03,BILW05,BIW10}.
Here we only deal with the case of $\PU(2,1)$. Among those representations, an important family is those admitting equivariant holomorphic maps.
These representations underly weight two variations of Hodge structures, exist provided the Toledo invariant is non-positive, and have been parametrized in Loftin--McIntosh~\cite{LM19b}.
We lack of a criterion, or a description in general of those holomorphic maps which are embeddings.
In this regard, a sufficient condition can be written in terms of the second fundamental form.
This criterion comes originally from Uhlenbeck~\cite{Uhl83}, and has been considered for Lagrangian immersions by
Loftin--McIntosh~\cite{LM13},
and can also be applied to holomorphic maps (see~\cite{Bro23Hadamard}). Namely, if the second fundamental form of a complete holomorphic immersion is bounded by $\eta<1$,
then it is properly embedded in the complex hyperbolic space. A representation admitting such an equivariant immersion
will then always be convex-cocompact, so in particular discrete and faithful. We call these representations \emph{almost-fuchsian}.
Our main result is the existence of almost-fuchsian representations admitting holomorphic equivariant immersions with nonmaximal Toledo invariant:
\begin{theorem*}[A]
	Let $d>0$ and $\eta>0$. There is a genus $g_0$ such that for every closed surface of genus $g>g_0$, there
	exists a representation $\rho:\pi_1(\Sigma_g)\rightarrow\PU(2,1)$ almost-fuchsian, admitting
	an equivariant holomorphic map $f$, verifying:
	\begin{equation*}
		\sup\|\II_f\|^2<\eta\quad\text{ and }\quad \rm{Tol}(\rho)=2-2g+\frac{2}{3}d\,.
	\end{equation*}
\end{theorem*}
This yields examples of nonmaximal representations which are nonetheless convex-cocompact.
Goldman--Kapovich--Leeb~\cite{GKL01} have constructed examples of convex-cocompact representations with any integer Toledo invariant, but this is
the first example of convex-cocompact representation in $\PU(2,1)$ which have non-integer Toledo invariant, hence don't lift to $\SU(2,1)$.
\begin{corollary}
	Provided the genus $g$ is large enough, there are convex-cocompact representations of a genus $g$ surface in $\PU(2,1)$ which
	do not lift to $\SU(2,1)$.
\end{corollary}
Another application is the construction of holomorphic immersions with small second fundamental form, yet which are not totally geodesic.
These holomorphic maps bound a quasicircle in $\partial_\infty\H^2_\C$, which is of Hausdorff dimension arbitrarily close to $1$, if $\eta$
is chosen close to zero.
\begin{corollary}
	For any $\eta>0$, there is a holomorphic embedding $\H^1_\C\rightarrow\H^2_\C$ whose second fundamental form verifies $\sup\|\II_f\|^2<\eta$,
	yet $f$ is not totally geodesic.
\end{corollary}
This corollary hints at the analogous question for holomorphic embeddings $\H^k_\C\rightarrow\H^{2k}_\C$.
Cao--Mok~\cite{CM90} have proven that any holomorphic immersion from $\H^n_\C$ to $\H^m_\C$ is totally geodesic if $m<2n$.
Koziarz--Maubon~\cite{KM17} have proven the rigidity of maximal representation for higher dimensional complex hyperbolic
lattices, yet the question remains open for the case of nonmaximal representations of lattices of $\PU(k,1)$, $k\geq 2$.
\subsection{Comparison with the case of $\SO(4,1)$}
In the paper~\cite{Bro23H4}, the author proved the existence of almost-fuchsian representations in $\SO(4,1)$,
whose corresponding hyperbolic $4$-manifold is diffeomorphic to the total space of a degree $1$ line bundle,
and which admits an equivariant superminimal immersion $f$ with small second fundamental form. While this shares some
similarities with the case studied here, the argument presented here
is more synthetic and does not rely on a study of the Moser--Trudinger inequality, enabling us to deal directly with
degree $d$ disc bundles rather than degree $1$. We also establish a sufficient non-asymptotic criterion for the existence
of almost-fuchsian equivariant immersions,
which was not discussed in the $\H^4$-case. The methods applied here actually also apply to $\H^4$, and show the existence
of almost-Fuchsian representations admitting superminimal equivariant embeddings, such that the uniformized hyperbolic
$4$-manifold is diffeomorphic to a degree $d$ disc bundle over the surface, as is discussed in the appendix.
\subsection{Some Context: Complex hyperbolic geometry and almost-Fuchsian representations}
\subsubsection{convex-cocompact representations in $\PU(2,1)$}
Xia~\cite{Xia00} has counted the connected components of the character variety of a surface group in $\PU(2,1)$.
He showed that they are indexed by their Toledo invariant, denoted $\rm{Tol}$, which belongs to $\frac{2}{3}\Z\cap[2-2g,2g-2]$.
Among these representations, a special family is that of convex-cocompact representations. These are discrete and faithful representations, leaving invariant a convex subset
of $\H^2_\C$ on which they act cocompactly. These representations all admit equivariant minimal surfaces.
In~\cite{GKL01}, Goldman--Kapovich--Leeb have shown that for all components with Toledo invariant in $\Z\cap[2g-2,2g-2]$, there are convex-cocompact representations with this prescribed Toledo
invariant. They also make a construction for noninteger Toledo invariant, but it doesn't lead to faithful representations.
Note that an integer Toledo invariant corresponds to when the representation lifts to $\SU(2,1)$. This lead Loftin--McIntosh~\cite{LM19b} to ask the question
\begin{question}
	Are there convex-cocompact representations valued in $\PU(2,1)$ which do not lift to $\SU(2,1)$ ?
\end{question}
We bring a positive answer to that question, at least when the genus of the surface is large enough. It remains unclear yet if every component of the character variety
contains convex-cocompact representations. Note that from results of Tholozan--Toulisse~\cite{TT21}, the analogous question is not true for the punctured sphere case.
\subsubsection{almost-Fuchsian representations}
Almost-Fuchsian representations were first studied by Uhlenbeck~\cite{Uhl83}. It is there defined as a representation from a surface group into $\PSL(2,\C)$
admitting an equivariant minimal surface in $\H^3$ whose principal curvatures are uniformly in $(-1,1)$. She then proceeded to show that such a representation is always
kleinian, i.e. discrete and faithful, and that the equivariant minimal surface is unique.
Since then, almost-fuchsian representations have had a rich history,\cite{Eps86,GHW10,Sep16,ES22}, with applications to several important geometric problems, such
as the counting of minimal surfaces~\cite{KM12,KW21,CMN22,Jia22,LN24} or to the foliation problem of hyperbolic $3$-manifolds~\cite{CMS23}.
A natural question we are interested in is the following:
\begin{question}
	Which topological manifolds can be obtained as $\rho\backslash X$ for $\rho:\pi_1 S\rightarrow \Isom(X)$ an almost-fuchsian representation ?
\end{question}
It turns out, when $X$ is a rank $1$ symmetric space, that such a manifold is always diffeomorphic to the total space of a vector bundle over the surface $S$.
When $X$ is $\H^3$, the $3$-manifold is always $S\times\R$. However, in~\cite{Bro23H4}, we showed that if the surface $S$ has genus large enough,
the degree $1$ disc bundle over $S$ can be obtained as the quotient of $\H^4$ by an almost-Fuchsian representation. Here we will show a slightly more general statement in $\PU(2,1)$,
obtaining degree $1-g+d$ disc bundles over $S$, provided $g$ is large enough.
\subsubsection{holomorphic equivariant immersions in $\PU(2,1)$}
From the non-abelian Hodge correspondence, there is an identification betweeen equivariant branched minimal surfaces in $\H^2_C$ and polystable $\PU(2,1)$-Higgs bundles whose Higgs field $\varphi$
satisfies $\Tr(\varphi^2)=0$. Hence, there is a $\C^\ast$-action on that moduli space defined by $z\cdot(\calE,\varphi)=(\calE,z\varphi)$.
Simpson~\cite{Sim91} characterized the fixed points of that action as \emph{complex variations of Hodge structures}. In the case of $\PU(2,1)$, they come in several flavors, described
by Loftin--McIntosh~\cite{LM19b}. As these representations are assumed to have special geometric significance, it is natural to look at almost-Fuchsian representations in them,
in the same way as we did in \cite{Bro23H4}.
Loftin--McIntosh~\cite{LM13} already looked at Lagrangian almost-Fuchsian representations, corresponding to weight $3$ variations of Hodge structures, and these are all Toledo invariant $0$ 
representations and are deformation of discrete and faithful representations in $\PO(2,1)\hookrightarrow\PU(2,1)$.
A weight $2$ variation of Hodge structure would correspond to a holomorphic almost-Fuchsian representation, that is a representation admit a holomorphic or anti-holomorphic equivariant minimal surface.
From the results of Toledo~\cite{Tol79}, we cannot deform representations in $\PU(1,1)\hookrightarrow\PU(2,1)$ to obtain irreducible representations. So we have to look at nonmaximal
representations, and these cannot be deformations of representations admitting a totally geodesic equivariant immersion.
Note that every connected component of the character variety contains a weight $2$ or a weight $3$ variation of Hodge structure, since the Energy functional is proper, see for instance Bradlow--Garcia-Prada--Gothen~\cite{BGG03}, Prop. 4.23.
\subsection{Acknowledgments}
The author is thankful for the support of the Max Planck Institut for Mathematics in the Sciences (MPIMIS) in his research. He is also thankful to Nicolas Tholozan and Andrea Seppi
for their help and insights over this topic.
This project has received funding from the European Research Council (ERC) under the European Union's Horizon 2020 research and innovation programme (grant agreement No 101018839).

\section{Holomorphic immersions in the complex hyperbolic plane}
\subsection{The curvature equations}
Let $S$ be a connected Riemann surface, not necessarily compact, and consider $f:\widetilde S\rightarrow \H^2_\C$
be a holomorphic immersion equivariant for some representation $\rho:\pi_1 S\rightarrow\Isom(\H^2_\C)=\PU(2,1)$.
Then the first and second fundamental form satisfy some relations, which we could call the "holomorphic Gauss--Codazzi equations".
We present here an analytic way of understanding holomorphic curves in $\H^2_\C$.
Note first that it is a classical result that the induced metric on a holomorphic submanifold will be negatively curved, hence $\widetilde S\approx \H^2$ is biholomorphic to the Poincaré disc. We will denote by $g_{h}$ its hyperbolic metric.
Let $L$ be a holomorphic line bundle over $S$. Consider $g_L$ a Hermitian metric on $L$ and $\lambda\in\R$ some constant
such that $F(g_L)=\lambda (-i\omega)$, where $\omega$ is the hyperbolic volume form.
\begin{theorem}\label{thm:fundaholoimm}
	Let $\beta\in H^0(S,K^3L^{-1})$, and $u,v:S\rightarrow\R$ be smooth functions satisfying
\begin{equation}
\left\{\begin{array}{ll}
	\Delta u&= 2e^{2u}-1+e^{-4u}e^{-2v}\|\beta\|_{hyp}^2\\
	\Delta v&= \lambda-3e^{2u}
\end{array}\right.
\end{equation}
Then there is a holomorphic immersion $f:\widetilde S\rightarrow\H^2_\C$, equivariant for some representation $\pi_1S\rightarrow\PU(2,1)$
satisfying:
\begin{enumerate}
	\item The induced metric is a lift of $e^{2u}g_{hyp}$.
	\item The normal bundle to $f$ is diffeomorphic to the universal cover of the total
		space of $KL^{-1}$, and its induced metric is $(e^{2u}g_{hyp})^\ast\otimes(e^{2v}g_L)^\ast$
	\item The $(2,0)$-part of the second fundamental form of $f$ is a lift of $\beta$.
\end{enumerate}
Conversely, for any holomorphic equivariant immersion $\widetilde S\rightarrow\H^2_\C$, these equations are satisfied
for the induced metrics on the tangent and normal vector bundles.
\end{theorem}
\begin{proof}
Assume that $(u,v,\beta)$ satisfy the aforementioned properties.
Denote by $\omega$ the hyperbolic volume form.
Consider the vector bundle $E=K^{-1}\oplus KL^{-1}\oplus\calO$,
endowed with the metric $H=e^{2u}g_{hyp}\oplus (e^{2u}g_{hyp})^\ast)(e^{2v}g_L)^\ast\oplus 1$,
$1$ denoting the standard flat metric on $\calO$.
The the following connection on $E$ is projectively flat:
\begin{equation}
D=\left(\begin{array}{ccc}
		\nabla & -\beta^\ast & 1\\
		\beta & \nabla & 0\\
		1^\ast & 0 &\nabla
		\end{array}\right)
\end{equation},
	where $\nabla$ is used to denote the Chern connection of the corresponding Hermitian line bundle.
	Note that here $\beta\in H^0(K^3L^{-1})$ is considered as a holomorphic $(1,0)$-form valued in $K^2L^{-1}$,
	in the same way that $1\in H^0(\calO)$ is thought of as a holomorphic $(1,0)$-form valued in $K^{-1}$.
Indeed, we just develop the following identities:
\begin{flalign}
F(K^{-1},e^{2u}g_{hyp})&=(-\Delta u -1)(-i\omega)\,,\\
F(L,e^{2v}g_L)&=(-\Delta v+\lambda)(-i\omega)\,,\\
\beta\wedge\beta^\ast&=e^{-4u}e^{-2v}\|\beta\|_{hyp}^2(-i\omega)\,,\\
	1\wedge 1^\ast&= e^{2u}(-i\omega)\,,
\end{flalign}
And then it comes:
\begin{equation}
	F(E,D)=\left(\begin{array}{lll} F(K^{-1})-\beta^\ast\beta+11^\ast&0&0\\
	0&\beta^\ast\beta+F(KL^{-1})&0\\
	0&0&F(\calO)+1^\ast 1\end{array}\right)=
	-11^\ast\otimes\one_{\End(E)}\,.
\end{equation}
Hence $(E,D)$ is a projectively flat bundle.
	Moreover, $(E,D)$ reduces to a $\PU(2,1)$-principal bundle, since the sesquilinear form of signature $(2,1)$
	given by $B$, admitting fiberwise the expression:
\begin{equation}
	B=\left(\begin{array}{lll}1&0&0\\0&1&0\\0&0&-1\end{array}\right)\,,
\end{equation}
is parallel, ($DB=0$).
We deduce that the developing map given by the metric is an equivariant map $f:\widetilde S\rightarrow\H^2_\C$,
and it verifies (up to lifting to the universal cover):
\begin{equation}
	f^\ast\nabla^{\H^2_\C}=\left(\begin{array}{lll}
		\nabla & -\beta^\ast & 0\\
		\beta & \nabla & 0\\
		0 & 0 & \nabla \end{array}\right)\, , (df)_\C=\left(\begin{array}{lll}0&0&1\\ 0&0&0\\ 1^\ast&0&0\end{array}\right)\,.
\end{equation}
The complex structure on $\H^2_\C$ induces the splitting of $T\H^2_\C=T^{1,0}\H^2_\C\oplus T^{0,1}\H^2_\C$,
with $f^\ast T^{1,0}\H^2_\C=\Hom(\calO,K^{-1}\oplus KL^{-1})$ and $f^\ast T^{0,1}\H^2_\C=\Hom(K^{-1}\oplus KL^{-1},\calO)$.
Here one can see that the $(1,0)$-part of $f$ is valued in $T^{1,0}\H^2_\C$, hence the map $f$ is a holomorphic immersion,
	equivariant with regard to the holonomy representation $\pi_1 S\rightarrow\PU(2,1)$.

The induced metric is the one pullback by $df$, so it is the one on the first factor $K^{-1}$ , hence a lift of $e^{2u}g_{h}$.

The complexification of the normal bundle is the intersection of the orthogonal of $T^\C\H^2_\C=\Hom(\calO,K^{-1})\oplus\Hom(K^{-1},\calO)$ with $T^\C\H^2_\C$, so it is $\Hom(\calO,KL^{-1})\oplus\Hom(KL^{-1},\calO)$.
Its $(1,0)$-part is $KL^{-1}$, endowed with the metric $(e^{2u}g_{hyp})^\ast\otimes(e^{2v}g_L)^\ast)$,	as claimed.

	Finally, the second fundamental form is the part of $f^\ast\nabla^{\H^2_\C}$ valued in the endomorphisms from
	the tangent bundle to the normal bundle, so here it is a lift of $\beta$, as desired.
\end{proof}
\begin{remark}[Comparison with the Higgs bundles parametrization]
\begin{enumerate}
\item
With the given notations, the Higgs bundle $(E,\dbar,\Phi)$ will satisfy the projective Yang-Mills equations,
	where $\Phi=1\in H^0(K\Hom(\calO,K^{-1}))$ and
\begin{equation}
\dbar=\left(\begin{array}{ccc}\dbar&-\beta^\ast&0\\0&\dbar&0\\0&0&\dbar\end{array}\right)\,.
\end{equation}
In that case, the Donaldson--Uhlenbeck--Yau proof of the non-abelian Hodge correspondence states that the existence of a solution to
the projective Yang--Mills equations is equivalent to the poly-stability of $(E,\dbar,\Phi)$. This leads to a parametrization of 
the space of equivariant holomorphic maps, and even of the equivariant minimal surfaces in $\H^2_\C$, as it was made in
Loftin--McIntosh~\cite{LM19b}. Here our parametrization is really in terms of germs of first and second fundamental form, and
we have no general result on existence. However, we will be able in some cases to explicitly bound the second fundamental form and get some
nice geometric properties.
\item
A priori, there is no reason to consider only the case when $u$ and $v$ are bounded. However, we will assume it to be the case, as
by results of Ahlfors and Yau, $u$ bounded corresponds to completeness of the induced metric, and when $v$ is assumed bounded, the
$\lambda$ such that $F(L,g_L)=\lambda(-i\omega)$ is uniquely defined.
Now in general, the geometric significance of the existence of $\lambda$ and
	$g_L$ such
that the conformal factor $e^{2v}$ is bounded is unknown.
\item
When $S$ is a closed Riemann surface of genus $g$, then $\lambda=\frac{\deg(L)}{2g-2}$
can take only discretely many values, and is related to the Toledo invariant of the holonomy
representation by the formula $\rm{Tol}=-\frac{2}{3}\deg(L)$. Also, we see that the presented
bundle is projectively equivalent to a flat $\SU(2,1)$-bundle if and only if $L$ is of degree divisible by $3$,
which is equivalent to the Toledo invariant being an integer.
\end{enumerate}
\end{remark}
\subsection{Stability conditions from the associated Higgs bundle}
Let $(E,\nabla,H)$ denote our projectively flat bundle with $H=h_{K^{-1}}\oplus h_{KL^{-1}}\oplus 1$ the hermitian metric involved.
Decomposing our connection as a sum of $H$-unitary and self-adjoint part,
we get:
\begin{equation}
	\nabla_E=\left(\begin{array}{ccc}
	\partial&0&0\\ \beta&\partial&0\\ 0&0&\partial\end{array}\right)+
		\left(\begin{array}{ccc}
	\dbar& -\beta^\ast &0\\ 0&\dbar&0\\ 0&0&\dbar\end{array}\right)+
	\left(\begin{array}{ccc}
	0&0&1\\0&0&0\\0&0&0\end{array}\right)+
		\left(\begin{array}{ccc}
		0&0&0\\ 0&0&0\\ 1^\ast &0&0\end{array}\right)
\end{equation}
Denoting this decomposition $\nabla=\partial+\dbar+\theta+\theta^\ast$,
$(E,\dbar,\theta)$ is the corresponding Higgs bundle to the harmonic metric $H$.
(The harmonicity of $H$ is encapsulated in $\dbar\theta=0$).
\begin{proposition}
	Assume that $\beta^\ast\neq 0$.
	If the Higgs bundle $(E,\dbar,\theta)$ is stable, then the degree of $L$ satisfies $0<\deg(L)< 3g-3$.
\end{proposition}
\begin{proof}
	Because of the shape of the Higgs field, invariant line subbundles are the ones
	in the kernel of the Higgs field, so subbundles
	of $K^{-1}\oplus KL^{-1}$. A particular holomorphic subbundle in the kernel is $K^{-1}$.
	It is a holomorphic subbundle of slope $2-2g$, while the slope
	of $E$ is $\frac{-1}{3}\deg(L)$.
	As such a first necessary condition for the slope-stability is:
	\begin{equation}
		\deg(L)\leq 6g-6\,.
	\end{equation}
	A rank $2$ subbundle $\theta$-invariant is either $K^{-1}\oplus KL^{-1}$,
	or the direct sum of a line subbundle of $K^{-1}\oplus KL^{-1}$ and of $\calO$.
	Among these, examples of holomorphic subbundles are $K^{-1}\oplus KL^{-1}$ and
	 $K^{-1}\oplus\calO$.
	For the Higgs bundle to be stable, it is then necessary that:
	\begin{equation}
		\left\{\begin{array}{rl}
			-\frac{1}{2}\deg(L)&<-\frac{1}{3}\deg(L)\\
			1-g&<-\frac{1}{3}\deg(L)
		\end{array}\right.
	\end{equation}
	Hence if it is stable, it satisfies $0<\deg(L)<3g-3$.
\end{proof}
\begin{remark}
	\begin{enumerate}
		\item
	When $\deg(L)=3g-3$, the representation is maximal, and $K^{-1}\oplus\calO$ is an invariant Higgs bundle of same slope. Hence the only possibility to obtain a polystable projective Higgs bundle is when $\beta=0$. This is symptomatic
	of the rigidity of maximal representations, which are all valued in
	a copy of $\PU(1,1)$ in $\PU(2,1)$. This statement is due to Toledo,\cite{Tol79}.
		\item
	When $\deg(L)=0$, there is no representation admitting an equivariant
			holomorphic immersion. However, there is another special
			family of immersions to consider, called \emph{Lagrangian immersions}. Almost-fuchsian Lagrangian immersions were already studied
			by Loftin--McIntosh~\cite{LM13}.
			It turns out their study is quite similar to
			the study of minimal surfaces in $\H^3$, as it boils
			down to the study of the Donaldson functional,
			see Huang--Lucia--Tarantello~\cite{HLT23}.
		\item
		To characterize stability rather than only having a necessary condition, one needs to 
		understand maximal subbundles of the rank $2$ bundle $V$ corresponding to the class $[\beta^*]$ and given
			by the exact sequence
			\begin{equation}
				0\rightarrow K^{-1}\rightarrow V\rightarrow KL^{-1}\rightarrow 0\,.
			\end{equation}
			A systematic study of these has been done by Lange and Narasimhan~\cite{LN83}, see also McIntosh~\cite{McI23}.
		In the construction we will present, we will prove by hand the existence of a projectively flat metric
		with the desired properties, which implies that the corresponding Higgs bundle is stable.
		\item
	With this convention, one can see that when $\beta=0$, the sectional curvature of a holomorphic totally geodesic plane in $\H^2_C$ is $-2$, so the sectional curvature of $\H^2_\C$ is pinched between $-2$ and $-\frac{1}{2}$.
	\end{enumerate}
\end{remark}

\section{Almost-Fuchsian immersions in the complex hyperbolic plane}
In this section we study a specific class of representations in $\PU(2,1)$,
called \emph{almost-Fuchsian}. Loftin and McIntosh already
considered the specific notion of Lagrangian almost-Fuchsian immersions,
but most of the story holds without considering the Lagrangian assumption.
Let $X$ denote a surface. Denote by $g_{\H^2_\C}$ the symmetric metric
on $\H^2_\C$, normalized so that the sectional curvature of $\H^2_C$
is bounded between $-4$ and $-1$.
\begin{definition}\label{def:AFimm}
Let $X$ be an orientable connected surface. An immersion $f:X\rightarrow\H^2_\C$
is said to be almost-Fuchsian if $f^\ast g_{\H^2_\C}$ is complete
and $\sup|\II_f|<1$.
\end{definition}
The main theorem is an application of \cite{Bro23Hadamard}.
\begin{theorem}
	Let $f:X\rightarrow\H^2_\C$ be an almost-Fuchsian immersion.
	Then $X$ is diffeomorphic to a disc,
	$(X,f^\ast g_{\H^2_\C})$ is quasi-isometric to the hyperbolic plane
	$\H$, $f$ is an embedding, a quasi-isometric
	embedding~$\H\rightarrow\H^2_\C$ and the exponential map of $f$
	defines a diffeomorphism
	\begin{equation}
		\exp_f:N_f X\approx\H^2_\C\,.
	\end{equation}
\end{theorem}
As a corollary, we get a sufficient criterion for convex-cocompactness
of a representation.
\begin{corollary}
	Let $S$ be a closed, orientable, connected surface, and $\rho:\pi_1 S\rightarrow\PU(2,1)$
	be a representation admitting an equivariant almost-Fuchsian immersion.
	Then $\rho$ is convex-cocompact.
\end{corollary}
Note that the completeness assumption on $f$ is redundant with the equivariance,
since $S$ is assumed compact.
The immersion assumption on $f$ cannot be weakened, as the examples of nonfuchsian representations
preserving a copy of $\PU(1,1)$ in $\PU(2,1)$. Those may admit a totally geodesic equivariant holomorphic map, but it will not be an immersion, and the representation cannot be convex-cocompact.

When that equivariant almost-Fuchsian immersion can be made minimal,
we say the representation is almost-Fuchsian:
\begin{definition}
	A representation $\rho:\pi_1 S\rightarrow\PU(2,1)$ is said to be \emph{almost-Fuchsian} if it admits an equivariant, minimal, almost-Fuchsian immersion
	$f:\widetilde S\rightarrow\H^2_\C$.
\end{definition}
Almost-Fuchsian representations admit a unique equivariant minimal surface.
This is a classical result from Uhlenbeck~\cite{Uhl83} for Kleinian representations,
see also El-Emam--Seppi~\cite{ES22}, here we deduce it from Theorem 1.3. of \cite{Bro23Hadamard}.
\begin{theorem}
	Let $\rho:\pi_1 S\rightarrow\PU(2,1)$ be an almost-Fuchsian representation.
	Denote by $f$ the equivariant minimal almost-Fuchsian immersion $f$.
	Then the image of $f$ is the unique $\rho(\pi_1 S)$-invariant minimal surface
	in $\H^2_\C$.
\end{theorem}
\subsection{Constraints on almost-fuchsian representations}
As far as we know, a characterization of the connected components containing almost-Fuchsian representations does not exist, neither does one of the connected components containing an almost-Fuchsian representation with holomorphic
equivariant map, which is the topic here.
To deal with almost-fuchsian representations, it will be useful to have a scalar interpretation of the curvature equations, requiring some notations, that will be used throughout this paper.
\begin{notation}
	Let $S$ denote a Riemann surface of genus $g>1$, with hyperbolic metric $h$.
	We denote by $\Delta$ its Laplace--Beltrami operator, and by $\omega$ the hyperbolic volume form.
	For any Hermitian line bundle $L$ over $S$ of nonzero degree, we denote by $h_L$ its Hermitian-Einstein metric,
	that is the unique metric whose curvature form satisfies:
	\begin{equation}
		\frac{i}{2\pi}F(L,h_L)=\frac{\deg L}{2g-2}\omega\,.
	\end{equation}
	For $\alpha$ a holomorphic section of $L$, we denote by $|\alpha|^2$ the square of its pointwise norm with regard to the metric $h_L$.
\end{notation}
The curvature equations of an almost-fuchsian holomorphic immersion can be re-written
as an additional boundedness assumption on our curvature equations:
\begin{theorem}\label{thm:fundaholoimmAF}
	Let $L$ be a Hermitian line bundle over $\Sigma$, $\beta\in H^0(S, K^3 L^{-1})$
	and $u,v:\Sigma_g\rightarrow\R$ be smooth solutions of 
	\begin{equation}
		\left\{\begin{array}{ll}
			\Delta u&=2e^{2u}-1+e^{-4u}e^{-2v}|\beta|^2\\
			\Delta v&=\frac{\deg(L)}{2g-2}-3e^{2u}
		\end{array}\right.
	\end{equation}
		with the additional control
	\begin{equation}
		\sup e^{-6u}e^{-2v}|\beta|^2\leq\eta<\frac{1}{2}\,.
	\end{equation}
		Then there is an almost-fuchsian
		representation $\rho:\pi_1\Sigma_g\rightarrow\PU(2,1)$
		with equivariant holomorphic map $f:\widetilde\Sigma_g\rightarrow\H^2_\C$
		whose $(2,0)$-part of the second fundamental form is a lift of $\beta$.
		The Toledo invariant of the representation $\rho$ depends on the degree of $L$
		in the following way:
	\begin{equation}
		\rm{Tol}(\rho)=-\frac{2}{3}\deg(L)\,.
	\end{equation}
\end{theorem}
\begin{remark}
	The constant $\frac{1}{2}$ is only because we work with an ambient metric whose sectional curvature is upper bounded
	by $\frac{-1}{2}$. If we had worked with the other normalization, we would have recovered the constant $1$ for
	the second fundamental form.
\end{remark}
\begin{proof}
	First we apply the fundamental theorem for holomorphic immersions, Theorem~\ref{thm:fundaholoimm}. This yields an equivariant holomorphic immersion $f:\widetilde\Sigma_g\rightarrow\H^2_\C$,
	whose $(2,0)$-part of the second fundamental form is a lift of $\beta$.
	Then by construction $e^{-6u}e^{-2v}|\beta|$ is the norm of the second fundamental form
	for the induced metric by $f$, hence it being less than  $\eta<1$ implies that
	$f$ is an almost-fuchsian immersion, and the representation $\rho$ for which
	it is equivariant is almost-fuchsian.

	It only remains to compute the Toledo invariant of $\rho$. This can be done directly
	from the expression of the projectively flat bundle, which as a $\PU(2,1)$-bundle
	splits as $\calU\oplus\calV$ with $\calU=K^{-1}\oplus KL^{-1}$ and $\calV=\calO$.
	Hence the standard expression of the Toledo invariant (see for instance~\cite{BGG03}). 
	\begin{equation*}
		\rm{Tol}(\rho)=\frac{2}{3}(\deg(\calU)-2\deg(\calV))=-\frac{2}{3}\deg(L)\,.
	\end{equation*}
\end{proof}

Almost-fuchsian representations with equivariant holomorphic maps can appear only when the Toledo invariant is close to
being maximal, as the following Proposition makes precise:
\begin{proposition}
	Let $\rho:\pi_1\Sigma\rightarrow\PU(2,1)$ be an almost-Fuchsian
	representation admitting an equivariant holomorphic map $f$.
	Then the Toledo invariant of $\rho$ satisfies
	\begin{equation*}
		2-2g\leq\rm{Tol}(\rho)<\frac{4}{5}(2-2g)\,.
	\end{equation*}
\end{proposition}
\begin{proof}
	Let $\rho:\pi_1\Sigma\rightarrow\PU(2,1)$ denote such a representation.
	Recall that $g>1$ is the genus of $\Sigma$.
	With the previously introduced notations, there are $u$,$v$ solutions
	of
	\begin{equation*}
		\left\{\begin{array}{cl}
		\Delta u &=2e^{2u}-1+e^{-4u}e^{-2v}|\beta|^2\\
		\Delta v &=\frac{\deg(L)}{2g-2}-3e^{2u}
		\end{array}\right.
	\end{equation*}
	Because it is almost-Fuchsian, we know that $e^{-6u}e^{-2v}|\beta|^2<\frac{1}{2}$.
	In particular, $u$ satisfies the following inequalities:
	\begin{equation*}
		2e^{2u}-1\leq \Delta u\leq \frac{5}{2}e^{2u}-1\,,
	\end{equation*}
	and by the maximum and minimum principles, we deduce that
	\begin{equation*}
		\frac{2}{5}\leq e^{2u}\leq\frac{1}{2}\,.
	\end{equation*}
	Now recall the equation on the curvature of $L$:
	\begin{equation*}
		\Delta v=\frac{\deg(L)}{2g-2}-3e^{2u}
	\end{equation*}
	Integrating it,
	we obtain
	\begin{equation*}
		\frac{6}{5}2\pi(2g-2)\leq 2\pi\deg(L)=3\int e^{2u}\omega\leq\frac{3}{2}2\pi(2g-2)\,.
	\end{equation*}
	Since $\rm{Tol}(\rho)=-\frac{2}{3}\deg(L)$, we obtain the corresponding bounds on the Toledo number:
	\begin{equation*}
		2-2g\leq\rm{Tol}(\rho)\leq\frac{4}{5}(2-2g)\,.
	\end{equation*}
	It only remains to check the strict inequality. If $\rm{Tol}(\rho)=\frac{4}{5}(2-2g)$,
	then it must be that $e^{2u}$ is constant equal to $\frac{5}{2}$. In particular, it follows that $\beta$ has
	no zeroes, which is absurd, hence the strict inequality.
\end{proof}

\section{The curvature equations of holomorphic immersions}
We present here some analysis of the curvature equations,
necessary for the existence theorems we will prove after.
\subsection{The Gauss Equation of a holomorphic immersion}
Let $f:\H\rightarrow\H^2_\C$ be a holomorphic immersion.
With the notations of the previous section, the curvature of the induced metric
satisfies an equation of the type:
\begin{equation}\label{GaussEq}
\Delta u =2e^{2u}-1+e^{-4u}f\,.
\end{equation}
where $f$ is a positive function. This equation can be solved when
$f$ is small enough, thanks to the sub- and supersolution method.
\begin{proposition}
Let $\alpha>0$, $\eta\in (0,1)$ and $f$ be an $\alpha$-Hölder function on $\H$
satisfying:
\begin{equation}
0\leq f \leq\frac{\eta}{(2+\eta)^3}\,.
\end{equation}
Then there is a unique $C^2$-regular solution $u$ to Equation \ref{GaussEq}
satisfying
\begin{equation}
e^{-6u}f\leq\eta\,.
\end{equation}
Furthermore, it satisfies the estimates
\begin{flalign}
\frac{-\ln(2+\eta)}{2}&\leq u\leq\frac{-\ln(2)}{2}\\
|\Delta u|_\infty&\leq\frac{\eta}{2+\eta}
\end{flalign}
\end{proposition}
\begin{proof}
The constant function $-\frac{\ln 2}{2}$ is a supersolution to equation \ref{GaussEq}.
Because $f\leq\frac{\eta}{(2+\eta)^3}$, one can check that the constant
function $\frac{-\ln(2+\eta)}{2}$ is a subsolution to the same equation.
	By the sub-supersolution method, there exists a solution $u$ to \ref{GaussEq}
satisfying
\begin{equation*}
\frac{-\ln(2+\eta)}{2}\leq u\leq\frac{-\ln(2)}{2}\,.
\end{equation*}
In particular, it verifies $e^{-6u}f\leq\eta$ and the other controls
claimed.

It remains to prove uniqueness of such a function.
If $u$ and $v$ are both $C^2$-regular solutions of the stated problem,
Then $w=u-v$ satisfies the inequality:
\begin{equation*}
\Delta w \geq 2e^{2v} w -fe^{-4v}w\,.
\end{equation*}
So in particular,
\begin{equation*}
\Delta w\geq \eta e^{2v} w\,.
\end{equation*}
	Now $w$ is assumed bounded, so applying the Omori--Yau maximum principle~\cite{Omo67,Yau75},
we deduce that
\begin{equation*}
w\geq 0
\end{equation*}
Interchange the roles of $u$ and $v$, and one gets $w\leq 0$.
Hence $w=0$ and $u=v$, as claimed.
\end{proof}
By the elliptic regularity principles, solving this equation
yields a continuous operator $F:U\subset C^{0,\alpha}(\H)\rightarrow C^{2,\alpha}(\H)$.
\begin{definition}
	Let $\calU$ be the following subset of $\C^{2,\alpha}(\H)\times C^{0,\alpha}(\H)$:
	\begin{equation}
		\calU=\{(u,f): f\geq 0,\,,\Delta u =2e^{2u}-1+e^{-4u}f,\sup e^{-6u}f<1\}\,.
	\end{equation}
\end{definition}
\begin{proposition}\label{Phiregular}
	The map $\Phi:\calU\rightarrow C^{0,\alpha}(\H)$, defined by $\Phi(u,f)=f$
	is a homeomorphism onto its image.
\end{proposition}
\begin{proof}
	We have already proven that $\Phi$ is injective. It only remains to check that $\Phi$ is a local homeomorphism.
	Let $\calF$ be the functional defined on the Banach space $C^{2,\alpha}(\H)\times C^{0,\alpha}(\H)$:
	\begin{equation*}
		\calF(u,f)=-\Delta u +2e^{2u}-1+e^{-4u}f\,.
	\end{equation*}
	Then its partial derivative has the explicit expression:
	\begin{equation*}
		\partial_u\calF(u,f)\dot u=-\Delta \dot u+4e^{2u}(1-e^{-6u}f)\dot u= L\dot u
	\end{equation*}
	We claim that, when $(u,f)$ belongs to $\calU$, $L$ is a linear isomorphism between $C^{2,\alpha}(\H)$ and $C^{0,\alpha}(\H)$.

	First, $L$ is injective: If $L\dot u=0$, then $\dot u$ is bounded and because $1-e^{-6u}f>0$, an application of Omori--Yau's maximum principle
	shows that $\dot u=0$.

	$L$ is also surjective. Let $w\in C^{0,\alpha}(\H)$. Choose $z_0\in\H$ and let $B_r(z_0)$ denote the open ball of hyperbolic radius $r$ around $z_0$.
	Because $1-e^{-6u}f>0$, the problem $Lv=w$ can be solved in any compact set, cf Gilbarg--Trudinger~\cite{GT01}, Theorem 6.8..
	Hence for every $r>0$ there exists $v_r\in C^{2,\alpha}_0(B_r(z_0))$ such that $Lv_r=w|_{B_r(z_0)}$.
	Applying the Schauder estimates to $L$ enables us to state the following:
	There is $C>0$ such that, for any $s>r+1$,
	\begin{equation*}
		\|v_s\|_{C^{2,\alpha}(B_r(z_0))}\leq C(\|v_s\|_{C^0}+\|L v_s\|_{C^{0,\alpha}})
	\end{equation*}
	Again, from the Omori--Yau maximum principle we show that
	\begin{equation*}
		\|v_s\|_{C^0}\leq\frac{2+\eta}{4\eta}\|Lv_s\|_{C^0}
	\end{equation*}
	Hence we deduce that, provided $s>r+1$:
	\begin{equation*}
		\|v_s\|_{C^{2,\alpha}(B_r(z_0))}\leq C(1+\frac{2+\eta}{4\eta})\|w\|_{C^0,\alpha}\,.
	\end{equation*}
	Applying the Arzela--Ascoli compactness principle, we get a subsequence of $v_r$ which converges on every compact set in the $C^{2,\beta}$-topology ($\beta<\alpha$)
	towards some limit $v\in C^{2,\alpha}(\H)$. This limit function then necessarily verifies
	\begin{equation*}
		Lv=w\,.
	\end{equation*}
	So we have proven that $L$ is bijective.
	By the closed graph theorem, $L$ is then a linear isomorphism, hence $\partial_u\calF$ is invertible. Applying the implicit function theorem,
	this directly shows that $\Phi$ is a local homeomorphism, as desired.
\end{proof}
\subsubsection*{Study of a ray of solutions}
Here we consider a ray of solutions $(u_t)$ to the PDEs
\begin{equation}
	\Delta u_t=2e^{2u_t}-1+e^{-4u_t}tf\,.
\end{equation}
We go on with a concavity statement concerning the volume of the metrics $e^{2u_t}g_h$. This is heavily inspired from
the work in \cite{BS24}, in which the authors parametrize the space of almost-fuchsian discs by a convex set of quadratic differentials.
Here we consider a ray of solutions, and show a concavity property of the volume.
While this can be deduced from \cite{BS24}, we give another proof here, as the context is slightly different.
\begin{proposition}
	Assume that for every $t>0$, $e^{-6u_t}tf< 1$, and that $f$ is not the zero function.
	Then the function
	\begin{equation}
		F(t)=\int_\Sigma e^{2u_t}\omega\,,
	\end{equation}
	is a nonincreasing concave function on $[0,1]$.
\end{proposition}
\begin{proof}
	The fact that $(u_t)$ is a smooth path in $C^{2,\alpha}(\Sigma)$
	can be deduced from the local invertibility of the equation and the uniqueness of solutions with the sup norm less than one.
	Differentiating in $t$, we get the control
	\begin{flalign}
		\Delta\dot u_t&=4(e^{2u_t}-e^{-4u_t}t f)\dot u_t+e^{-4u_t}f\\
		\Delta\ddot u_t&=4(e^{2u_t}-e^{-4u_t}t f)\ddot u_t+4(2e^{2u_t}+4e^{-4u_t}tf)\dot u_t^2-8 t e^{-4u_t}f \dot u_t
	\end{flalign}
	Now the control on the first derivative implies
	\begin{equation*}
		\Delta\dot u_t\geq 4e^{2u_t}(1-|e^{6u_t}tf|_\infty)\dot u_t\,,
	\end{equation*}
	Hence by the Omori--Yau maximum principle, $\dot u_t\leq 0$.
	Now this means that the second derivative satisfies
	\begin{equation*}
		\Delta\ddot u_t\geq 4e^{2u_t}(1-|e^{-6u_t}tf|_\infty)\ddot u_t\,,
	\end{equation*}
	and so by the maximum principle, $\ddot u_t\leq 0$. But this won't be
	a sufficient control to obtain concavity of the volume.
	Using the convexity property of the Laplacian,
	we get the inequality
	\begin{flalign*}
		\Delta(\ddot u_t+2\dot u_t^2)&\geq\Delta\ddot u_t+4\dot u_t\Delta\dot u_t\\
		&\geq 4(e^{2u_t}-e^{-4u_t}tf)\ddot u_t+24 e^{2u_t}\dot u_t^2
		-8te^{-4u_t}f\dot u_t+e^{-4u_t}f\\
		&\geq 4(e^{2u_t}-e^{-4u_t}tf)(\ddot u_t+2\dot u_t^2)+16 e^{2u_t}\dot u_t^2
	\end{flalign*}
	From the maximum principle we deduce that $\ddot u_t+2\dot u_t^2\leq 0$.
	Finally, we must write the second derivative of the volume:
	\begin{equation*}
		F''(t)=\int_\Sigma 2e^{2u_t}(2\dot u_t^2+\ddot u_t)d\omega
	\end{equation*}
	Hence because $2\dot u_t^2+\ddot u_t\leq 0$, $F$ is concave, as desired.
	
	Finally, we can check that $F'(0)$ has the following expression:
	\begin{equation*}
		F'(0)=\int_\Sigma 2\dot u_0 e^{2u_0}d\omega=\int_\Sigma \dot u_0 d\omega\,.
	\end{equation*}
	We have already shown that $\dot u_0\leq 0$, and it is straightforward that $\dot u_0=0$ if and only if $f=0$.
	Hence whenever $f$ is nonzero, $F'(0)<0$ and $F$ is nonincreasing, as desired.
\end{proof}
This concavity gives us a precious bound on the volume of solutions
to Gauss equations:
\begin{proposition}
	Assume that for every $t$, $|e^{-6u_t}t f|_\infty<1$.
	Then the volume $F(1)$ satisfies:
	\begin{equation}
		\fint e^{2u_1}d\omega\leq\frac{1}{2}-2t\fint fd\omega\,.
	\end{equation}
\end{proposition}
\begin{proof}
	We have already proven that $F$ is concave. Hence applying the slope inequality, we get
	\begin{equation*}
		F(1)\leq F(0)+tF'(0)\,.
	\end{equation*}
	Now $u_0$ is the constant function such that $e^{2u_0}=\frac{1}{2}$,
	hence $F(0)=\frac{\Vol(\Sigma)}{2}$.
	Also, $\dot u_0$ is solution of
	\begin{equation*}
		\Delta\dot u_0=2\dot u_0+4f\,.
	\end{equation*}
	Hence the following value of $F'(0)$
	\begin{equation*}
		F'(0)=\int_\Sigma \dot u_0d\omega=-2\int_\Sigma fd\omega\,.
	\end{equation*}
	Finally, dividing by the volume of the surface we obtain
	\begin{equation*}
		\fint e^{2u_1}d\omega\leq\frac{1}{2}-2t\fint fd\omega\,,
	\end{equation*}
	as claimed.
\end{proof}
This leads naturally to the consideration of the following ratio.
\begin{definition}[Balance ratio]
	Let $(\Sigma,h)$ be a closed hyperbolic surface, compatible with the Riemann Surface Structure $S$.
	Let $f:\Sigma\rightarrow\R_+$ be a measurable, bounded function.
	Then we call the Balance ratio of $f$ the following quantity:
	\begin{equation}
		\rm{bal}(f)=\frac{\fint f}{|f|_\infty}\,.
	\end{equation}
	Let $L\rightarrow S$
	be a holomorphic line bundle over $\Sigma$
	and $\alpha\in H^0(L)$. Call $h_L$ the uniformizing Hermitian-Einstein metric on $L$.
	Then we call the Balance ratio of $\alpha$ the quantity:
	\begin{equation}
		\rm{bal}(\alpha)=\frac{\fint|\alpha|_{h_L}^2}{|\alpha|^2_{\infty,h_L}}\,.
	\end{equation}
\end{definition}
A direct corollary of our discussion is that the Gauss equation can be solved with prescribed volume provided the given data is sufficiently balanced:
\begin{corollary}\label{cor:GaussVolume}
	Let $R>0$ and $\eta\in(0,1)$. Consider $f\in C^{0,\alpha}(\Sigma)$
	a nonnegative function such that:
	\begin{equation}
		\rm{bal}(f)\geq\frac{(2+\eta)^3 R}{2\eta}\,.
	\end{equation}
	Then there exists a unique $u\in C^{2,\alpha}(\Sigma)$ and $t>0$
	such that
	\begin{flalign}
		\Delta u &= 2e^{2u}-1+e^{-4u}tf\\
		\sup e^{-6u}tf&\leq\eta\\
		\fint e^{2u}d\omega&=\frac{1}{2}-R
	\end{flalign}
\end{corollary}
This corollary shows the interest into finding balanced sections of Hermitian line bundles. Note that the analogous statement for the Gauss equation
of a minimal surface in $\H^3$, $\Delta u =e^{2u}-1+e^{-2u}f$, will show
that a control on the balance ratio allows to consider almost-Fuchsian surfaces
with defect $R$ in the volume of the induced metric. This ratio also appeared
in \cite{Bro23H4} in the construction of almost-Fuchsian structures on
degree $1$ disc bundles over a surface. Here we will use it to build examples
of almost-Fuchsian representations in $\PU(2,1)$ with nonmaximal Toledo invariant.
\begin{lemma}[Regularity of the solution]\label{Psiregular}
	Let $R>0$ and $\eta\in(0,1)$ as in Corollary \ref{cor:GaussVolume}.
	Let $U_{R,\eta}$ denote the subset in $C^{0,\alpha}(\Sigma)$:
	\begin{equation}
		U_{R,\eta}=\{ f\in C^{0,\alpha}(\Sigma):f\geq 0,\,\rm{bal}(|f|)>
			\frac{(2+\eta)^3 R}{2\eta}\}\,.
	\end{equation}
	Then the induced map
	\begin{equation}
		\Psi:U_{R,\eta}\rightarrow C^{2,\alpha}(\Sigma)\times\R
	\end{equation}
	which to $f$ gives the solution $(u,t)$ from Corollary \ref{cor:GaussVolume} is continuous.
\end{lemma}
\begin{proof}
	Let $f\in U_{R,\eta}$ and $(u,t)=\Psi(f)$.
	Because $f$ is nonzero, we know that the volume function is strictly nonincreasing along rays,
	and it is a submersion at $f$, and this shows that the second projection of $\Psi$
	$t=p_2\circ\Psi(f)$ is continuous on $U_{R,\eta}$.
	But $u=p_1\circ\Psi(f)$ can be obtained as $p_1\circ\Phi^{-1}\circ (p_2\circ\Psi(f)\cdot f)$, with $\Phi$
	the functional considered in Proposition~\ref{Phiregular}, which we have proven to be a local homeomorphism.
	Hence we conclude that $\Psi$ is continuous, as desired.
\end{proof}
\subsection{Study of the Poisson Equation}
We don't pretend here to give an exhaustive study of the Poisson equation on a hyperbolic surface.
We just state the results we will use, with some idea of the proofs, which are very classical.
Recall $(\Sigma,h)$ be a closed hyperbolic surface. We consider $f$ a positive function,
$r>0$ and we wish to solve the equation
\begin{equation}\label{PoissonEq}
	\Delta v = r -f\,.
\end{equation}
Contrasting with the study of the Gauss Equation, the estimates we will give here will depend more on the intrinsical geometry of the surface.
In particular, two quantities will be extensively used:
\begin{notation}[Systole, Spectral Gap]
	We denote by $\delta$ the systole of $(\Sigma,h)$, that is the length of the shortest nontrivial closed curve in $(\Sigma,h)$.
	We denote by $\Lambda$ the spectral gap of $(\Sigma,h)$, that is the largest constant such that, for any zero-average function $u\in C^{1}(\Sigma)$:
	\begin{equation}
		\|u\|_2^2\leq\Lambda^{-1}\|\nabla u\|^2_2\,.
	\end{equation}
\end{notation}
A necessary condition for the existence of a solution,
is that
\begin{equation}\label{VolEq}
	\int_\Sigma fd\rm{Vol}(h)=r\Vol(\Sigma,h)\,.
\end{equation}
\begin{proposition}\label{prop:Poissonsolve}
	Let $f\in C^{0,\alpha}(\Sigma)$ and $r>0$ satisfying
	condition~\ref{VolEq}. Then there is a unique zero-average function 
	$v\in C^{2,\alpha}(\Sigma)$ satisfying
	\begin{equation}
		\Delta v = r -f\,.
	\end{equation}
	Moreover, its sup norm may be controlled in the following way:
	\begin{equation}
		|v|_\infty\leq C(\delta,\Lambda)|r-f|_2\,.
	\end{equation}
\end{proposition}
\begin{proof}
	The existence statement can easily proven by sub-and supersolution,
	because $f$ is bounded. The uniqueness statement is because any bounded 
	harmonic function on $\Sigma$ is constant, hence if it is zero-average, 
	it must vanish. We go on with the estimate of $|v|_\infty$.
	First, we use the Morrey--Sobolev embedding~$W^{2,2}(\Sigma)\rightarrow L^\infty(\Sigma)$. Because this statement is a local one and $\Sigma$ is
	a hyperbolic surface, this constant may be controlled only by
	the systole of the surface, $\delta$.
	It remains to show that
	\begin{equation*}
		|v|_{2,2}\leq C(\Lambda)|r-f|_2\,.
	\end{equation*}
	A short proof of that goes in the following way:
	By the Poincaré inequality,
	\begin{equation*}
		|v|_2^2\leq \Lambda^{-1} |\nabla v|_2^2\,.
	\end{equation*}
	Also, by the Bochner identity,
	\begin{equation*}
		|\nabla^2 v|_2^2=|\nabla v|_2^2+|\Delta v|_2^2\,.
	\end{equation*}
	All in all, this means
	\begin{equation*}
		|v|_{2,2}^2\leq |r-f|_2^2+(2+\Lambda^{-1})|\nabla v|_2^2\,.
	\end{equation*}
	Finally, the Cauchy-Schwarz inequality ensures that
	\begin{equation*}
		|v|_2^2\leq\Lambda^{-1}|\nabla v|_2^2\leq \Lambda^{-1} |v|_2|r-f|_2\,.
	\end{equation*}
	From which we deduce that $|v|_2\leq \Lambda^{-1}|r-f|_2$,
	and in turn
	\begin{equation*}
		|\nabla v|_2^2\leq\Lambda^{-2}|r-f|_2^2\,.
	\end{equation*}
	This shows the explicit control:
	\begin{equation*}
		|v|_{2,2}\leq \sqrt{1+2\Lambda^{-2}+\Lambda^{-3}}|r-f|_2\,.
	\end{equation*}
	Hence we get
	\begin{equation*}
		|v|_\infty\leq C(\delta,\Lambda)|r-f|_2\,,
	\end{equation*}
	as claimed.
\end{proof}
Once again, it is straightforward that the induced map $C^{0,\alpha}_0(\Sigma)\rightarrow C^{2,\alpha}(\Sigma)$ is continuous (smooth).
Conveniently, the image by this map of any bounded set is bounded in $C^{2,\alpha}(\Sigma)$,
which sits compactly in $C^{0,\alpha}(\Sigma)$, by Arzela--Ascoli.

\section{Construction of balanced sections of line bundles}
In this section, we consider the notion of balanced sequence of sections of line bundles:
\begin{definition}
	A sequence $(\Sigma_g,N_g,\alpha_g)_{g\geq 2}$,
	where $\Sigma_g$ is a genus $g$ hyperbolic surface,
$N_g$ is a Hermitian line bundle over $\Sigma_g$ and $\alpha_g\in H^0(N_g)$ is said to be \emph{balanced}
if the following conditions are satisfied:
\begin{flalign*}
	\inf\delta(\Sigma_g)&>0\\
	\inf\Lambda(\Sigma_g)&>0\\
	\inf\rm{bal}(\alpha_g)&>0\,.
\end{flalign*}
\end{definition}
To prove our main result we need the following existence statement:
\begin{theorem}\label{thm:balanced}
Let $d>0$, $n>0$
There exists a balanced sequence $(\Sigma_g,N_g,\alpha_g)$ such that
$N_g$ is of degree $n(g-1)+d$.
\end{theorem}
\begin{proof}
The construction is made in the following way:
Consider $\Sigma_2$ a genus $2$ hyperbolic surface, and $\Sigma_g\rightarrow\Sigma_2$ a degree $g-1$ cover of~$\Sigma_2$,
such that $\Lambda(\Sigma_g)>\eps$ for some genus-independent constant $\eps>0$.
Such a sequence exists as a corollary of the works of Magee--Naud--Puder~\cite{MNP22}.
By construction, its systole is always larger than the systole of $\Sigma_2$.

It only remains to construct a balanced family of line bundles and holomorphic sections of it,
which is the done in the following way:
	Denote by $K_\sigma$ the cotangent bundle of $\Sigma_g$, of degree $2g-2$.
	Let $s_2$ be a nonzero holomorphic section of $K_2^\frac{n}{2}$.
	If $n>1$, such a section always exists in virtue of Riemann--Roch's Theorem,
	if $n=1$ it might not exist for every choice of squareroot of the cotangent bundle,
	but we can always choose that squareroot so that its space of holomorphic sections
	is odd-dimensional, hence non-zero, cf Atiyah~\cite{Ati71}.
	Now, lift this section to a sequence $s_g$ of sections of $K_g^\frac{n}{2}$.
Finally, consider a sequence of points $z_g\in\Sigma_g$, and consider the line bundles
\begin{equation*}
	N_g=K^{\frac{n}{2}}_\sigma\calO(z_g)^d\,.
\end{equation*}
	We observe first that $N_g$ is of degree $n(g-1)+d$, as desired.
	Also, denoting by $f_g$ the defining holomorphic section of $\calO(z_g)$,
we get that $s_g f_g^d$ is a holomorphic section of $N_g$.
	It remains to prove that this is a balanced sequence.
This relies on the elliptic study of $|f_g|^2$, cf \cite{Bro23H4}, Proposition 3.9,
	that we report here for the reader's convenience:
\begin{proposition}\label{prop:fgbal}
Let $r>0$ be smaller than half the systole $\delta$ of $\Sigma$. Then there are constants $C_1,C_2$ depending only
	on $(\delta,\Lambda,r)$ and $\lambda>0$ such that
\begin{equation}
\left\{\begin{array}{cl}
|\lambda f_g|^2\leq C_2(r_g,\delta_g,\Lambda_g)&\text{ on }D(z_g,r)\\
\frac{1}{C_1(r_g,\delta_g,\Lambda_g)}\leq|\lambda f_g|^2\leq C_1(r_g,\delta_g,\Lambda_g)&\text{ on }\Sigma_g- D(z_g,r)
\end{array}\right.
\end{equation}
\end{proposition}
As a consequence, we can show that the sequence $(s_g f_g^d)$ is balanced.
Indeed, fix $r$ small enough so that every disk of radius $r$ in $\Sigma_g$ is the lift
of a radius $r$ disk in $\Sigma_2$.
Then, on one side, we can estimate the sup norm:
\begin{equation*}
\frac{1}{C_1^d}|s_g|^2_\infty\leq |s_g f_g^d|_\infty^2\leq \max(C_1,C_2)^d|s_g|^2_\infty\,.
\end{equation*}
And also we have the lower bound on the $L^2$-norm:
\begin{equation*}
\int_{\Sigma_g} |s_g f_g^d|^2 d\Vol(h)\geq \frac{(g-2)}{C_1^d}\int_{\Sigma_2}|s_2|^2\,.
\end{equation*}
In particular, this shows
\begin{equation*}
\rm{bal}(s_g f_g^d)\geq \frac{1}{C_1^{d}\max(C_1,C_2)^d}\frac{g-2}{g-1}\rm{bal}(s_2)\,.
\end{equation*}
Hence the sequence is balanced, as claimed.
\end{proof}

\section{Existence of almost-Fuchsian holomorphic maps}
In this section, we prove the following result:
\begin{theorem}\label{thm:exoAFholo}
	Let $d>0$, $\eta\in(0,\frac{1}{2})$. There is a genus~$g_0>1$ such
	that for any $g\geq \sigma_0$,
	there exists a representation~$\rho:\pi_1\Sigma_g\rightarrow\PU(2,1)$
	and an equivariant holomorphic
	map~$f:\widetilde\Sigma_g\rightarrow\H^2_\C$ satisfying:
	\begin{enumerate}
		\item
			The second fundamental form of $f$ satisfies $|\II_f|\leq\eta$.
		\item
			The Toledo invariant of $\rho$ is $2-2g+\frac{2d}{3}$.
	\end{enumerate}
\end{theorem}
As a consequence of this theorem, the representation $\rho$ is almost-Fuchsian, hence convex-cocompact,
hence discrete and faithful, yet lifts to $\SU(2,1)$ if and only if $d$ is a multiple of $3$.
Note that since $\rho$ admits a holomorphic equivariant immersion, its Toledo invariant must be non-positive. An analogous
result holds for representations admitting anti-holomorphic equivariant immersions, with non-negative Toledo invariant.
\begin{corollary}
	In sufficiently large genus, there are convex-cocompact representations
	of a genus~$g$ surface in $\PU(2,1)$ which do not lift to $\SU(2,1)$.
\end{corollary}
This answers a question raised by Loftin and McIntosh in \cite{LM13}.
The main ingredient of the proof of theorem \ref{thm:exoAFholo} is the following fixed-point theorem.
Let $(\Sigma,g)$ be a finite volume hyperbolic surface with systole $\delta$ and spectral gap $\Lambda$.
Denote by $C(\delta,\Lambda)$ the constant such that, for any zero-average function $v\in C^{2,\alpha}(\Sigma)$:
\begin{equation}
	|v|_\infty\leq C(\delta,\Lambda)|\Delta v|_2
\end{equation}
Let $\eta\in(0,\frac{1}{2})$ and $R>0$. Then we have the following criterion:
\begin{theorem}\label{thm:fixedpoint}
	Assume there exists $A>0$ and $f\in C^{0,\alpha}(\Sigma)$ a nonnegative function satisfying
	\begin{flalign}
		\rm{bal}(f)&\geq A\\
		A\cdot\exp\big(\frac{-12 C(\delta,\Lambda)\eta\sqrt{\Vol(\Sigma)}}{2(2+\eta)}\big)&\geq\frac{(2+\eta)^3 R}{2\eta}
	\end{flalign}
	Then there exists $u,v\in C^{2,\alpha}(\Sigma)$ and $t\in\R_+$ satisfying:
	\begin{equation}
		\left\{\begin{array}{l}
		\Delta u =2e^{2u}-1+e^{-4u}e^{-2v}tf\\
		\Delta v =\frac{3}{2}-3R-3e^{2u}\\
		\sup|e^{-6u}e^{-2v}tf|\leq\eta\,.
		\end{array}\right.
	\end{equation}
\end{theorem}
\begin{proof}
	Consider the following chain of continuous maps:
	\begin{equation*}
		U_1\overset{\Phi_1}{\rightarrow} U_2\overset{\Phi_2}{\rightarrow} U_3\,.
	\end{equation*}
	Where
	\begin{flalign*}
		U_1&=\{\hat f\in C^{0,\alpha}(\Sigma),\,\hat f\geq 0,\,\rm{bal}(\hat f)\geq\frac{(2+\eta)^3 R}{2\eta}\}\\
		U_2&=\{\hat u\in C^{0,\alpha}(\Sigma),\,\fint e^{2\hat u}=\frac{1}{2}-R,\,\frac{1}{2+\eta}\leq e^{2\hat u}\leq\frac{1}{2}\}\\
		U_3&=\{\hat v\in C^{0,\alpha}(\Sigma),\,\fint v=0,|v|_\infty\leq 3C(\delta,\Lambda)\frac{\eta\sqrt{\Vol(\Sigma)}}{2(2+\eta)}\}
	\end{flalign*}
	$U_1$ is a closed convex subset of $C^{0,\alpha}(\Sigma)$.
	$\Phi_1$ is the continuous map obtained from Lemma~\ref{Psiregular} which to such an $f$ associates $\hat u\in U_2$ satisfying
	\begin{equation*}
		\Delta \hat u =2e^{2\hat u}-1+e^{-4\hat u}\hat f\,.
	\end{equation*}
	$\Phi_2$ is the continuous map which to $\hat u$ associates the zero-average solution of the Poisson equation:
	\begin{equation*}
		\Delta\hat v=\frac{3}{2}-3R-3e^{2\hat u}\,.
	\end{equation*}
	we now exhibit a continuous map $\Phi_3:U_3\rightarrow U_1$.
	This map is the following: $\Phi_3(\hat v)=\hat f=e^{-2\hat v}f$.
	We need to estimate the Balance ratio of $\hat f$ to show that $\Phi_3(\hat v)$ belongs to $U_1$:
	\begin{equation*}
		\rm{bal}(\hat f)\geq e^{-4|\hat v|_\infty}\rm{bal}(f)\geq A\cdot\exp\big(\frac{-12 C(\delta,\Lambda) \eta\sqrt{\Vol(\Sigma)}}{2(2+\eta)}\big)
	\end{equation*}
	The main hypothesis of the theorem then ensures that
	\begin{equation*}
		\rm{bal}(\hat f)\geq \frac{(2+\eta)^3 R}{2\eta}
	\end{equation*}
	So $\hat f\in U_1$ as claimed.
	Also, due to elliptic regularity, the image of $\Phi_2$ has compact closure in $C^{0,\alpha}(\Sigma)$. Hence $\Phi_3\circ\Phi_2\circ\Phi_1$ has a fixed point in $U_1$,
	by the Banach--Schauder fixed point theorem.
	Denote $\tilde f$ this fixed point, and $u=\Phi_1(\tilde f)$, $v=\Phi_2\circ\Phi_1(\tilde f)$.
	Then the fixed point property means that
	\begin{equation*}
		\tilde f=e^{-2v}f\,.
	\end{equation*}
	And by construction, $u$ and $v$ satisfy, for some $t>0$:
	\begin{equation*}
		\left\{\begin{array}{cl}
			\Delta u&=2e^{2u}-1+e^{-4u}e^{-2v}tf\\
			\Delta v&= \frac{3}{2}-3R-3e^{2u}
		\end{array}\right.
	\end{equation*}
	Finally, by construction of the map $\Phi_1$, we have the desired upper bound:
	\begin{equation*}
		\sup|e^{-6u}e^{-2v}tf|\leq\eta\,.
	\end{equation*}
\end{proof}
Combined with the existence of nicely balanced sections shown in the previous sections, we are now equipped to prove Theorem~\ref{thm:exoAFholo}:
\begin{proof}[Proof of Theorem~\ref{thm:exoAFholo}]
	Fix $d>0$. From Theorem~\ref{thm:balanced},
	there exists a balanced sequence~$(\Sigma_g,N_g,\alpha_g)$ with~$N_g$
	a line bundle of degree $3g-3+d$ over~$\Sigma_g$,
	and $\alpha_g\in H^0(N_g)$.
	Denote by~$K_g$ the cotangent bundle
	to $\Sigma_g$, and by $L_g=K^{3}N_g^{-1}$ the line bundle of
	degree~$3g-3-d$.
	Denote by $\delta_0,\Lambda_0 $ and $A_0$ the nonzero constants such that, for every $g$,
	\begin{flalign*}
			\delta(\Sigma_g)&\geq\delta_0\\
			\Lambda(\Sigma_g)&\geq\Lambda_0\\
			\rm{bal}(\alpha_g)&\geq A_0\,.
	\end{flalign*}
	We want to apply Theorem~\ref{thm:fixedpoint} with $f=|\alpha_\sigma|^2$,
	and with the volume constant fixed by the degree of $L_\sigma$:
	It must satisfy:
	\begin{equation*}
		\frac{3}{2}-3R=\frac{\deg(L)}{2g-2}=\frac{3}{2}-\frac{d}{2g-2}\,.
	\end{equation*}
	We need to obtain $R=R_g=\frac{d}{6g-6}$.
	In order to do so, we need to find a sequence~$\eta_g\in(0,1)$,
	such that the condition
	\begin{equation*}
		A_0\cdot\exp\big(\frac{-12 C(\delta_0,\Lambda_0)\eta_g
		\sqrt{\Vol(\Sigma_g)}}{\sqrt{2(2+\eta_g)}}\big)
		\geq\frac{(2+\eta_g)^3 R_g}{2\eta_g}
	\end{equation*}
	will be verified, at least in the limit of~$g$ large.

	Because $\Vol(\Sigma_g)=2\pi(2g-2)$, we observe that if we pick a
	sequence~$\eta_g\in(0,1)$ such that:
	\begin{equation*}
		g\eta_g\rightarrow+\infty\quad\text{and}\quad
		\eta_g\sqrt{g}\rightarrow 0\,,
	\end{equation*}
	then the condition above will be verified, as the left-hand side converges to $A_0$
	while the right-hand side converges to $0$.
	Hence there is $g_0$ such that in genus larger than $g_0$ the above condition
	will be verified,
	and by Theorem~\ref{thm:fixedpoint} we get the existence of $t>0$ and smooth
	functions $u,v:\Sigma\rightarrow\R$ verifying:
	\begin{equation*}
		\left\{\begin{array}{rl}
			\Delta u &=2e^{2u}-1+e^{-4u}e^{-2v}t|\alpha_g|^2\\
			\Delta v &=\frac{3}{2}-\frac{d}{2g-2}-3e^{2u}\\
			\sup e^{-6u}t|\alpha_g|^2&\leq\eta_g<1
		\end{array}\right.
	\end{equation*}
		In particular, applying Theorem~\ref{thm:fundaholoimmAF},
		we get the existence of an almost-fuchsian
		representation~$\rho:\pi_1\Sigma_g\rightarrow\PU(2,1)$
		with an equivariant holomorphic embedding~$f:\widetilde\Sigma_g\rightarrow\H^2_\C$
		satisfying $\sup|\II_f|\leq\eta$,
		and whose Toledo invariant verifies $\rm{Tol}(\rho)=-\frac{2}{3}\deg(L_g)=2-2g+\frac{2d}{3}$,
		as desired.
\end{proof}
It is tempting to conjecture that this proof could be carried for a sequence of line bundles $(L_g)$
whose degrees satisfy the following asymptotics:
\begin{equation*}
	\frac{d_g}{\sqrt{g}}\rightarrow 0\,.
\end{equation*}
However, to complete that proof one would need the existence of balanced sections of $L_g$ with that
prescribed degree, which remains an open problem up to now.

\appendix
\section{Superminimal surfaces in $\H^4$}
In the paper~\cite{Bro23H4}, the author studied almost-Fuchsian representations
in $\SO_0(4,1)$. In particular, the following theorem was shown:
\begin{theorem}
	Let $\eta\in (0,1)$. There is a genus $g_0>0$ such that,
	for every genus $g>g_0$, there exists a representation
	$\rho:\pi_1\Sigma_g\rightarrow\SO_0(4,1)$ satisfying:
	\begin{enumerate}
		\item
			$\rho$ is almost-Fuchsian,
			with equivariant minimal map $f:\H^2\rightarrow\H^4$.
		\item
			The embedding $f$ is superminimal, i.e. its normal Gauss
			map is a conformal map to the space of geodesic discs of $\H^4$.
		\item
			$f$ satisfies $\|\II_f\|\leq\eta$.
		\item
			The hyperbolic manifold $\rho\backslash\H^4$
			is diffeomorphic to the total space of
			a degree $1$ disc bundle over $\Sigma_g$.
	\end{enumerate}
\end{theorem}
	While the argument in~\cite{Bro23H4} relies heavily on the degree $1$
	particularity, as one can then use the Moser--Trudinger to solve
	the Ricci equation, with the study carried out in this paper we can
	get rid of this degree 1 specificity, and show:
\begin{theorem}\label{thm:exoAFH4}
	Let $\eta>0$ and  $d>0$. There is a genus $g_0>0$, such that
	for every $g>g_0$, there exists a representation $\rho:\pi_1\Sigma_g\rightarrow\SO_0(4,1)$ satisfying:
	\begin{enumerate}
		\item
			$\rho$ is almost-Fuchsian, with equivariant minimal map
			$f$ satisfying $\|\II_f\|\leq\eta$.
		\item
			$f$ is a superminimal map in $\H^4$.
		\item
			The hyperbolic manifold $\rho\backslash\H^4$
			is diffeomorphic to the total space of a degree $d$
			disc bundle over $\Sigma_g$.
	\end{enumerate}
\end{theorem}
The idea is that we can prove this theorem in the same way as we constructed
nonmaximal almost-Fuchsian representations with equivariant holomorphic maps.
Indeed, to get an equivariant superminimal immersion in $\H^4$
we need to find a solution of the following PDE system (see Proposition 2.11. of \cite{Bro23H4})
\begin{proposition}\label{prop:fundthmAFH4}
	Let $\eta>0$
	Let $\Sigma$ be a closed hyperbolic surface of genus $g$. Denote by $\omega$
	its volume form. Consider $N$ a line bundle of degree $d\geq 0$ endowed with its uniformizing
	metric $h_N$ of curvature form $c\omega=\frac{d}{2g-2}\omega$. Let $\alpha\in H^0(K^2N)$
	and $u,v$ smooth functions on $\Sigma$ satisfying:
	\begin{equation}
		\left\{\begin{array}{l}\label{eq:H4uv}
			\Delta u= e^{2u}-1+e^{-2u}e^{2v}|\alpha|^2\\
			\Delta v= c- e^{-2u}e^{2v}|\alpha|^2\\
			\sup e^{-4u}e^{2v}|\alpha|^2\leq\eta<1\,.
		\end{array}\right.
	\end{equation}
		Then there is a convex-cocompact representation $\rho:\pi_1\Sigma\rightarrow\SO_0(4,1)$
		and a minimal, superminimal, $\rho$-equivariant and $\eta$-almost-fuchsian
		immersion $f$ such that:
	\begin{enumerate}
		\item
			The induced metric by $f$ is a lift of $e^{2u}\omega$.
		\item
			The hyperbolic manifold $\rho\backslash\H^4$ is diffeomorphic to
			a degree $d$ disc bundle over $\Sigma$.
		\item
			The induced metric on the normal bundle to $f$ is $e^{2v}h_N$.
		\item
			The holomorphic second fundamental form of $f$ is a lift of $\alpha$.
		\end{enumerate}
\end{proposition}
In order to prove Theorem~\ref{thm:exoAFH4}, we need to show the existence, in
genus large enough, of $u,v$ solutions to Equations~\ref{eq:H4uv}. The trick is then
to introduce $w=u+v$. Hence we see that $u,v$ is a solution to Equations~\ref{eq:H4uv}
if and only if $u,w$ is a solution to these Equations:
\begin{equation}\label{eqH4uw}
	\left\{\begin{array}{l}
		\Delta u=e^{2u}-1+e^{-4u}e^{2w}|\alpha|^2\\
		\Delta w=c-1+e^{2u}\\
		\sup e^{-6u}e^{2w}|\alpha|^2\leq\eta<1\,.
	\end{array}\right.
\end{equation}
We then mimic the proof of the existence of almost-Fuchsian equation.
On the first equation, we can write the following Theorem
\begin{lemma}\label{lemma:GaussH4map}
	Let $R>0$ and $\eta\in(0,\frac{1}{2})$.
	Let $f\in C^{0,\alpha}(\Sigma)$ belong to the set $V_{R,\eta}$:
	\begin{equation}
		V_{R,\eta}=\big\{f\in C^{0,\alpha}(\Sigma),\,f\geq 0\,
		\rm{bal}(f)\geq\frac{(8+\eta)^3}{16\eta}R\big\}\,.
	\end{equation}
	Then there exists a unique $u\in C^{2,\alpha}(\Sigma)$ and $t>0$ such that
	\begin{equation}
		\left\{\begin{array}{l}
			\Delta u = e^{2u}-1+e^{-4u}tf\\
			\sup e^{-6u}tf\leq\eta<1\\
			\fint e^{2u}=1-R\,.
		\end{array}\right.
	\end{equation}
	The resulting map $\Phi:V_{R,\eta}\rightarrow C^{2,\alpha}(\Sigma)$
	is continuous.
\end{lemma}
\begin{proof}
	This is a direct application of Corollary~\ref{cor:GaussVolume} 
	with data $\widetilde R=\frac{R}{2}$, $\widetilde\eta=\frac{\eta}{2}$
	and $\widetilde f=\frac{f}{4}$.
	We get then $\widetilde t>0$ and a function $\widetilde u$ satisfying:
	\begin{equation*}
		\left\{\begin{array}{l}
			\Delta \widetilde u=2e^{2\widetilde u}-1+e^{-4\widetilde u}(\frac{tf}{4})\\
			\sup e^{-6\widetilde u}t\widetilde f\leq\widetilde\eta\\
			\fint e^{2\widetilde u}=\frac{1}{2}-\widetilde R\,.
		\end{array}\right.
	\end{equation*}
	Then it is clear that $u=\widetilde u+\frac{\ln 2}{2}$ and $t=\widetilde t$
	satisfy the prescribed conditions.
	
	The regularity of the map $\Phi$ corresponds to the regularity statement
	of Lemma~\ref{Psiregular}.
\end{proof}
Also, remark that every element in the image of $\Phi$ satisfies
\begin{equation}
\frac{4}{4+\eta}\leq e^{2u}\leq 1\,.
\end{equation}
Then, adapting Theorem~\ref{thm:fixedpoint} to the $\H^4$-setup, we have the following criterion:
\begin{theorem}\label{thm:fixedpointH4}
	Let $0<\eta<\frac{1}{2}$, and $R>0$.
	Let $C(\delta,\Lambda)$ be the constant defined in Prop.~\ref{prop:Poissonsolve}
	Assume there exists $A>0$ and $f\in C^{0,\alpha}(\Sigma)$
	a nonnegative function satisfying:
	\begin{equation}
		\left\{\begin{array}{rl}
			\rm{bal}(f)&\geq A\\
		A\cdot\exp\big(-\frac{4C(\delta,\Lambda)\sqrt{\Vol(\Sigma)}\eta}{4+\eta}\big)
			&\geq \frac{(8+\eta)^3}{16\eta}R\,.
		\end{array}\right.
	\end{equation}
	Then there exist $(u,w)\in C^{2,\alpha}(\Sigma)$ and $t>0$ satisfying:
	\begin{equation}
		\left\{\begin{array}{l}
			\Delta u =e^{2u}-1+e^{-4u}e^{2w}tf\\
			\Delta w =R-1+e^{2u}\\
			\sup{e^{-6u}e^{2v}tf}\leq\eta<1\,.
		\end{array}\right.
	\end{equation}
\end{theorem}
\begin{proof}
We use a fixed-point argument, à la Banach-Schauder. Let $V_{R,\eta}$ and $\Phi=(\Phi_u,\Phi_t)$
be the functional constructed in Lemma~\ref{lemma:GaussH4map}.
The image of $\Phi_u$ is valued in the set:
\begin{equation*}
W_\eta=\{ u\in C^{2,\alpha}(\Sigma):\frac{4}{4+\eta}\leq e^{2u}\leq 1\}
\end{equation*}
In particular, we can bound the $L^2$-norm:
\begin{equation*}
|R-1+e^{2u}|_2\leq\frac{\eta}{4+\eta}\sqrt{\Vol(\Sigma)}\,.
\end{equation*}
Applying Prop.~\ref{prop:Poissonsolve}, we know there is a unique zero-average function $w=\Psi(u)$
in $C^{2,\alpha}(\Sigma)$
such that
\begin{equation*}
	\Delta w=r-1+e^{2u}\,,
\end{equation*}
and it satisfies
\begin{equation*}
|w|_\infty\leq C(\delta,\Lambda)\frac{\eta}{4+\eta}\sqrt{\Vol(\Sigma)}\,.
\end{equation*}
The resulting map $\Psi:W_\eta\rightarrow C^{2,\alpha}(\Sigma)$ is smooth, and
has bounded image in $C^{2,\alpha}(\Sigma)$, so relatively compact in $C^{0,\alpha}(\Sigma)$.

Finally, consider $F(w)=e^{2w}f$. The written condition ensures that $F\circ\Psi(W_\eta)\subset V_{R,\eta}$.
Hence $F\circ\Psi\circ\Phi_u$ is a continuous self map of the closed convex set $V_{R,\eta}$,
with relatively compact image. By the Banach--Schauder fixed point theorem,
it admits a fixed point $\hat f$. Denote by $u=\Phi_u(\hat f)$, $t=\Phi_t(\hat f)$
and $w=\Psi\circ\Phi_u(\hat f)$, then from the fixed point property we get
\begin{equation*}
\left\{\begin{array}{l}
	\Delta u = e^{2u}-1+e^{-4u}e^{2w}tf\\
	\Delta w = R-1+e^{2u}\\
	\sup{e^{-6u}e^{2w}tf}\leq\eta<1\,,
	\end{array}\right.
\end{equation*}
as desired.
\end{proof}
It remains to prove Theorem~\ref{thm:exoAFH4}
\begin{proof}[Proof of Theorem~\ref{thm:exoAFH4}]
Fix $d>0$, and $\eta>0$. Without loss of generality we may assume $\eta<\frac{1}{2}$.
Let $(\Sigma_g,L_g,\alpha_g)$ be a balanced sequence with $\deg(L_g)=4g-4-d$.
Denote $N_g=K_g^{-2}L_g$, so that $\alpha_g\in H^0(K^2 N_g)$.
Denote by $\delta_0,\Lambda_0$ the lower bounds on the systoles and spectral gaps of $\Sigma_g$.
Also, consider $A_0>0$ the infimum of the Balance ratios $\rm{bal}(\alpha_g)$.

We want to apply Theorem~\ref{thm:fixedpointH4} to $f_g=|\alpha_g|^2$
with $R_g=\frac{d}{2g-2}$.
To do so, we need to find $0<\eta_g<\eta$ such that:
\begin{equation*}
	A_0\cdot\exp\big(-\frac{4C(\delta_0,\Lambda_0)\sqrt{\Vol(\Sigma_g)}\eta_g}{4+\eta_g}\big)
	\geq\frac{(8+\eta_g)^3}{16\eta_g}\frac{d}{2g-2}\,.
\end{equation*}
The fact that $\sqrt{\Vol(\Sigma_g)}$ grows like $\sqrt{g}$ shows that if we choose $(\eta_g)$
a sequence satisfying:
\begin{equation*}
	\eta_g\cdot g\rightarrow +\infty\quad\eta_g\sqrt{g}\rightarrow 0\,,
\end{equation*}
then that condition will be satisfied for every $g\geq g_0$ large enough.
Without loss of generality, we can take $g_0$ large enough such that $\eta_g\leq\eta$ also.
Hence applying,Theorem~\ref{thm:fixedpointH4}, for every $g>g_0$ we get $u,w,t$ solutions
of
\begin{equation*}
	\left\{\begin{array}{l}
	\Delta u = e^{2u}-1+e^{-4u}e^{2w}t|\alpha_g|^2\\
	\Delta w = \frac{d}{2g-2}-1+e^{2u}\\
	\sup (e^{-6u}e^{2w}tf)\leq\eta_g\leq\eta\,.
	\end{array}\right.
\end{equation*}
Denoting $v=w-u$, we get that $(u,v,t)$ is a solution of
\begin{equation*}
	\left\{\begin{array}{l}
	\Delta u = e^{2u}-1+e^{-2u}e^{2v}t|\alpha_g|^2\\
	\Delta v = \frac{d}{2g-2}-e^{-2u}e^{2v}t|\alpha_g|^2\\
	\sup (e^{-4u}e^{2v}t|\alpha_g|^2)\leq\eta_g\leq\eta\,.
	\end{array}\right.
\end{equation*}
From Prop.\ref{prop:fundthmAFH4} we then get that there is $\rho:\pi_1\Sigma_g\rightarrow\SO_0(4,1)$
an $\eta$-almost-Fuchsian representation with equivariant superminimal immersion $f$,
and whose normal bundle is a degree $d$ disc bundle, proving our Theorem.
\end{proof}

\newpage
\bibliographystyle{alpha}
\bibliography{references}
\end{document}